\documentclass [12pt] {article}
\usepackage{fullpage}
\usepackage{amssymb}
\usepackage{amsthm}
\usepackage{amsmath}
\usepackage{latexsym}
\usepackage{epsfig}
\usepackage{amscd}
\usepackage{graphicx}
\usepackage[dvips]{color}
\newtheorem{theorem}{Theorem}[section]
\newtheorem{lemma}[theorem]{Lemma}
\newtheorem{corollary}[theorem]{Corollary}
\theoremstyle{definition}
\newtheorem{definition}[theorem]{Definition}
\newtheorem{question}[theorem]{Question}

\newtheorem{remark}[theorem]{Remark}
\numberwithin{equation}{section}
\newcommand{\confdim}{\mathcal{C}\dim}
\newcommand{\eps}{\varepsilon}
\newcommand{\G}{\Gamma}
\newcommand{\g}{\gamma}
\renewcommand{\d}{\delta}
\newcommand{\la}{\lambda}

\newcommand{\bR}{\mathbb{R}}
\newcommand{\diam}{\mathrm{diam}}
\newcommand{\dist}{\mathrm{dist}}
\newcommand{\dmod}{\mathrm{d\textbf{-}mod}}

\newcommand{\calE}{\mathcal{E}}
\newcommand{\calB}{\mathcal{B}}

\renewcommand{\mod}{\mathrm{mod}}
\begin{document}

\title{Conformal dimension:\\ Cantor sets and moduli}
\author{Hrant Hakobyan}
\maketitle
\begin{abstract}
In this paper we give several conditions for a space to be minimal for conformal dimension.
We show that there are sets of zero length and conformal dimension $1$ thus answering a question of Bishop and Tyson.
Another sufficient condition for minimality is given in terms of a modulus of a system of measures in the sense of Fuglede ~\cite{Fug}. 
\begin{equation}
r^{1+\eps}\leq \lambda(E\cap B(x,r))
\end{equation}
It implies in particular that there are many sets $E\subset\mathbb{R}$ of zero length  such that $X\times Y$ is minimal for conformal dimension for every compact $Y$.
\end{abstract}


\section{Introduction}

Given a homeomorphism $\eta:[0,\infty)\rightarrow[0,\infty)$
a map $f$ between metric spaces $(X,d_X)$ and $(Y,d_Y)$ is called $\eta$-\textit{quasisymmetric} if for all distinct triples $x,y,z\in{X}$ and $t>0$
\begin{equation}\label{QS}
\frac{d_X(x,y)}{d_X(y,z)}\leq{t} \quad \Rightarrow \quad
\frac{d_Y(f(x),f(y))}{d_Y(f(y),f(z))}\leq{\eta(t)}.
\end{equation}
If $\eta(t)\leq C\max\{t^K,t^{1/K}\}$ for some $K\geq1$ and $C>0$ then $f$ is said to be \textit{power quasisymmetric}. We will denote by $\mathcal{QS}(X)$ the collection of all quasisymmetric maps defined on $X$.

 \textit{Conformal dimension} of a metric space, a concept introduced by Pansu in ~\cite{Pansu}, is the infimal Hausdorff dimension of quasisymmetric images of $X$,
$$\confdim{X}=\inf_{f\in\mathcal{QS}(X)}\dim_H f(X).$$ We say $X$ is \textit{minimal for conformal dimension} if
$\confdim{X}=\dim_H{X}$. Euclidean spaces with standard metric are the simplest examples of  minimal spaces. The first
examples of minimal sets of non integer dimension $\geq1$ were given in ~\cite{Pansu2} and ~\cite{Brd}. The minimality
in this examples was due to the presence of certain families of curves. In ~\cite{Tyson} Tyson proved that if $X$ is an Ahlfors $Q$-regular space then $\confdim{X}\geq Q$ if there is a curve family $\G$ in $X$ of positive $Q$ modulus (see Section \ref{sectionModulus} for the definitions and the statement of Tyson's theorem). In particular
$(0,1)\times Y$ is minimal for every Borel metric space $Y$.
The first minimal Cantor sets were constructed in ~\cite{BT01}. These Cantor sets were of Hausdorff dimension
$\geq1$ and had infinite Hausdorff $1$-measure. On the other hand in ~\cite{Kovalev} Kovalev proved a conjecture of
Tyson that if $\dim_H X<1$ then $\confdim{X}=0$. In ~\cite{Hak} the author proved that middle interval Cantor sets are minimal if one
considers quasisymmetric maps of the line to itself. In ~\cite{HuWen} this result was generalized to include a
larger class of uniform Cantor sets (see Section $2$  for the definitions). One of the main results of this paper, Theorem \ref{Thm1}, gives a sufficient condition for a metric space to have conformal dimension at least 1. The following Theorem is a consequence of Theorem \ref{Thm1}, see Remark \ref{RemA}, and answers a question of Bishop and Tyson from ~\cite{BT01}.
\begin{theorem}\label{A}
There is a set $E\subset\bR$ of zero length and conformal dimension $1$.
\end{theorem}
Theorem \ref{Thm1} also generalizes the main result of ~\cite{HuWen}, see Remark \ref{HuWen}.
In ~\cite{BT01} it was also shown that $E\times{Y}$ is minimal for every compact $Y\subset\bR^n$ if the Hausdorff $1$-contents of quasisymmetric images of $E$ are uniformly bounded away from $0$.
One of the main results of this paper, Theorem \ref{Thm2}, is a sufficient condition for a space $X$ to be minimal in terms of a certain modulus of a system of measures in the sense of Fuglede ~\cite{Fug}. It implies that given $E\subset\bR$ the products $E\times{Y}$ are minimal for compact $Y$ if $E$ is minimal and supports a measure with certain growth property. 

\begin{theorem}\label{B}
If  $E\subset\mathbb{R}$ is minimal and  supports a measure $\lambda$ s.t. for every $\eps>0$
$$r^{1+\eps}\lesssim \lambda(E\cap B_r(x))\lesssim r^{1-\eps}$$ for all $x\in E$ and all $r>0$  then  $E\times{Y}$ is minimal for every nonempty compact $Y$.
\end{theorem}

In the same article Bishop and Tyson asked for a characterization of subsets $E$ of the line which have the property that the product of $E$ with every compact $Y$ is minimal. It is clear that to have this property $E$ would have to be minimal itself. Theorem \ref{B} indicates that the converse may also be true. So the following is a natural question.
\begin{question}
Is $E\times{Y}$ minimal for every compact $Y$ if and only if $E\subset\bR$ is minimal?
\end{question}
Theorem \ref{B} does not quite answer this question since a lower bound on the Hausdorff dimension does not in general imply that there is a measure $\lambda$ satisfying the condition of the theorem, even though by Frostman's lemma there is a measure on $E$ which satisfies a growth estimate from above, namely for every $t<1$ and every ball of radius $r$  one has $\lambda(E\cap B_r)\lesssim r^t$.

Examples of sets $E$ which satisfy the conclusions of Theorems \ref{A} and \ref{B} are easy to construct. Consider the so called \textit{middle interval Cantor sets} constructed as follows. Start from the unit interval on the line. Remove its $c_1$-st middle part to obtain two intervals of equal length.
By induction, in the $i$-th step remove $c_i$-th middle part of every remaining
component from the previous step to obtain $2^i$ intervals of equal length.
If $c_i\to0$ and $\sum_{i\geq1}c_i=\infty$ then the resulting Cantor set
$E$ would satisfy the conclusions of Theorems \ref{A} and \ref{B}.
In fact we will show that all uniformly perfect middle interval Cantor sets are minimal if (and only if) they have Hausdorff dimension $1$.

In would  also be interesting to know whether it is necessary for one of the sets $X$ or $Y$ to be minimal in order for the product $X\times{Y}$ to be minimal. In view of Kovalev's Theorem
an easier question is the following. Are there two sets $X$ and $Y$ of  dimension $<1$ such that $\confdim{X\times Y}\geq1$?

This paper is organized as follows. In Section \ref{Background} we provide some background material and fix the notations. In section $3$ we state Theorem \ref{Thm1}  and explain how Theorem \ref{A} follows from it which . In Section \ref{sectionThm1proof} we proof Theorem \ref{Thm1}. In Section \ref{sectionModulus} we recall the definitions of the modulus of a system of measures and discrete modulus and deduce Theorem \ref{B} from Theorem \ref{Thm2}. We prove Theorem \ref{Thm2} in Section \ref{sectionThm2}.

\section{Background}\label{Background}
Constants in this article will be denoted by the letter $C$ and can have different values from line to line.
The notation $A\lesssim B$ means there is a constant $C$ such that $A\leq C B$. Given $r>0$ by $B_r$ we will denote any open ball in $X$ of radius $r$ and by $B(x,r)$ the one centered at $x\in X$ and by $CB(x,r)$ we will denote the ball $B(x,Cr)$.

Recall that the Hausdorff $t$-measure of a metric space $(X,d_X)$ is
defined as follows. For every open cover $\{U_i\}_{i=1}^{\infty}$ of $X$ let
\begin{equation*}
H^{\eps}_t(X)=\inf{ \left\{
\sum_{i=1}^{\infty}(\textrm{diam}{U_i})^t:
X\subset\bigcup_{i=1}^{\infty}U_i, \, \textrm{diam}{U_i}<\varepsilon
\right\} },
\end{equation*}
and $$H_t(X)=\lim_{\eps\to0} H^{\eps}_t(X).$$
 The
Hausdorff dimension of $X$ is
\begin{eqnarray*}
\dim_H(X)=\inf\{\,t\,: H_t(X)=0\}=\sup\{\,t\,: H_t(X)=\infty\}
\end{eqnarray*}

One usually gives an upper bound for the Hausdorff dimension of a
set by finding explicit covers for it. Lower bounds can be obtained
by finding a measure on $X$.

\begin{lemma}[Mass distribution principle]
If the metric space $(X,d_X)$ supports a positive Borel measure
$\mu$ satisfying $\mu(U)\leq{C(\diam{U})^{d}}$, for some fixed
constant $C>0$ and every $U\subset{X}$ then $\dim_H(E)\geq{d}$.
\end{lemma}
\begin{proof}For every cover $\{U_i\}_{i=1}^{\infty}$ of $X$ we have
$
\sum_i{(\diam{U_i})^d}\geq\frac{1}{C}\sum_i\mu(U_i)\geq\frac{1}{C}\mu(X).
$
Therefore $H^d(X)\geq\frac{\mu(X)}{C}>0$.
\end{proof}
An important converse is the following lemma, see ~\cite{Mattila}.

\begin{lemma}[Frostman's Lemma]
If $X$ is a metric space of Hausdorff dimension $d$ then there is a finite and positive measure $\mu$ on $X$ such that
$$\mu(B_r)\lesssim r^d.$$
\end{lemma}

\section{Conformal dimension of Cantor sets}
Let as first recall the following definition from \cite{SW}.

\begin{definition} \label{MIC} Given a sequence $\{c_i\}$ such that $0\leq c_i<1$, a set $E\subset\mathbb{R}$ is called \textit{$\{c_i\}$-thick} if there is a sequence of sets $\mathcal{E}_n=\{E_{n,j}\}$, where for each $n$ the $E_{n,j}$ are intervals with mutually disjoint interiors, such that $\sup_j\diam{E_{n,j}}\to0$, as $n\to\infty$ and each $E_{n,j}\setminus{E}$ contains an interval $J_{n,j}$ so that the following conditions are satisfied
\end{definition}
\begin{eqnarray*}
\frac{\diam{J_{n,j}}}{\diam E_{n,j}}&\leq& c_n, \\
\bigcup_{\calE_{n+1}}(E_{n,j}\setminus J_{n,j})&\subset& \bigcup_{\calE_{n+1}} E_{n+1,k},\\
\bigcup_{\calE_{n+1}}(E_{n+1,k}\setminus J_{n+1,k})&\subset& \bigcup_{\calE_n}(E_{n,j}\setminus J_{n,j}),\\
\bigcap_n\bigcup_j (E_{n,j}\setminus J_{n,j})&\subset& E.
\end{eqnarray*}

Note that a particular example of $\{c_i\}$ thick sets are the middle interval Cantor sets described in the introduction.

If $\sum_i c_i<\infty$ then a $\{c_i\}$-thick set has a positive Lebesgue measure on the line.
It was shown in \cite{SW} that if $\sum_i c_i^{p}<\infty$ for every $p>0$ then $E$ is quasisymmetrically thick, i.e. $f(E)$ has positive Lebesgue measure whenever $f:\bR\to\bR$ is a quasisymmetric map.
In the case of the middle interval Cantor sets the condition was shown to be necessary and sufficient for $E$ to be quasisymmetrically thick, see \cite{BHM}.

For every interval $E_{n,j}$ let $r_{n,j}$ denote the ratio of the lengths of the longer of the two components of $E_{n,j}\setminus J_{n,j}$ to the shorter one.

\begin{theorem}\label{Thm1}
Let $E$ be a $\{c_i\}$-thick set and $f$ a power quasisymmetric embedding of $E$ into some metric space. If
 \begin{eqnarray}
 \label{product} \sqrt[n]{\prod_{i=1}^n(1-c_i)}&\to&1, \mbox{and} \\
 \label{ratio} r_{n,j}&\leq& M, \mbox{  for some  } M<\infty
 \end{eqnarray}
then $\dim_H{f(E)}\geq1$.

\end{theorem}

\begin{corollary}\label{Cor1}
Suppose $E\subset\mathbb R$ is a middle interval Cantor set
 \begin{itemize}
 \item[\textit{(i).}]{If $E$ is uniformly perfect then it is minimal for conformal dimension if and only if $\dim_H E=1$.}
 \item[\textit{(ii).}]{If $\dim_H E=1$ then $\dim_H f(E)\geq1$ whenever $f$ extends to a quasisymmetric map of a uniformly perfect space.}
 \end{itemize}
\end{corollary}
 Recall that a metric space is \emph{uniformly perfect} if there is a constant
$C\geq1$ so that for each $x\in{X}$ and for all $r>0$
\begin{equation*}
{X\setminus{B(x,r)}\neq\emptyset\quad \Longrightarrow \quad
B(x,r)\setminus{B(x,\frac{r}{C})}\neq\emptyset}.
\end{equation*}
This condition in a sense rules out ``large gaps" in the space.
Examples of uniformly perfect sets are connected sets as well as
many totally disconnected sets, like middle third Cantor set or many
sets arising in conformal dynamics.
The importance of uniform perfectness in quasiconformal geometry comes from the following fact, see ~\cite{H}.
\begin{theorem}\label{ThmUP}
Any quasisymmetric embedding of a uniformly perfect space is
power-quasisymmetric.
\end{theorem}


\begin{proof}[Proof of Corollary \ref{Cor1}] By Kovalev's theorem for $(i)$ we only need to show that if $\dim_H E=1$ then $E$ is minimal.
Since every quasisymmetric map of a uniformly perfect space is power quasisymmetric and in the case of middle interval Cantor sets $r_{n,j}=1$ to prove $(i)$ and $(ii)$ we only need to show that $\dim_H E=1$ implies (\ref{product}).

Let $N(X,\varepsilon)$ be the minimal number of $\varepsilon$ balls
needed to cover $X$. Recall that \textit{upper and lower Minkowski dimensions} of $X$ are defined as
\begin{eqnarray*}
\overline{\dim}_M(X)=\limsup_{\varepsilon\rightarrow0}\frac{\log{N(X,\varepsilon)}}{\log{1/\varepsilon}} \,\,\,\,\mbox{  and  }\,\,\,\,
\underline{\dim}_M(X)=\liminf_{\varepsilon\rightarrow0}\frac{\log{N(X,\varepsilon)}}{\log{1/\varepsilon}} \label{dim}
\end{eqnarray*}
respectively. When these two numbers are the same the common value is
called Minkowski dimension of $X$ and is denoted by $\dim_M{X}$.
Generally  $\dim_H(X)\leq\underline{\dim}_M(X)\leq\overline{\dim}_M(X)$, see \cite{Mattila}.
Therefore if $X\subset\mathbb{R}$ and $\dim_H(X)=1$
then Minkowski dimension of $X$ exists, is equal $1$ and
\begin{eqnarray}
\dim_M(E)&=&\lim_{n\rightarrow\infty}\frac{\log{2^n}}{\log\frac{2^n}{\prod_{i=1}^n{(1-c_i)}}} \nonumber\\
&=&\lim_{n\rightarrow\infty}\frac{1}{1-\frac{1}{\log2}\log\sqrt[n]{\prod_{i=1}^n{(1-c_i)}}}=1.
\end{eqnarray}Therefore $\dim_H E=1$ if and only if (\ref{product}) holds.
\end{proof}

\begin{remark}\label{HuWen}
Theorem \ref{Thm1} generalizes the result of Hu and Wen from ~\cite{HuWen} where it was shown that $\dim_H f(E)=1$ whenever $E=E(\{n_i\},\{\g_i\})$ is a uniform Cantor sets of Hausdorff dimension $1$ corresponding to a bounded sequences $\{n_i\}$ and $f:\mathbb R\to\mathbb R$ is a quasisymmetric maps. Recall from ~\cite{HuWen} that given a sequence of positive integers $\{n_i\}$  and a sequence of real number $\{\g_i\}$ in $(0,1)$  a \textit{uniform Cantor set} $E$ corresponding to these sequences is constructed as follows. Divide $E_0=[0,1]$ into $n_1$ intervals of equal length so that the spacing between adjacent ``children" of $E_0$  is $\g_1\diam{E_0}$. In the $i$-th step divide every component $E_{i,j}$
remaining from the previous step into $n_i$ equal length intervals so that the distance between every two adjacent ones is $\g_i\diam E_{i,j}$.

It is not hard to see that $E$  satisfies the conditions of Theorem \ref{Thm1} if $n_i\leq N$ and $\dim_H E=1$. Therefore under these conditions if $E$ is uniformly perfect, which means $\g_i<C<1$, then $\confdim E=1$.
Also, even if $E$ is not uniformly perfect, $\dim_H f(E)\geq 1$ if $f$ extends to a quasisymmetry of a uniformly perfect space (for instance a quasiconformal map of a Euclidean space as in ~\cite{HuWen}).

The fact that $\{n_i\}$ is a bounded sequence is crucial in this case since otherwise one can easily construct a uniform Cantor set of Hausdorff dimension $1$ which does not satisfy the condition \ref{product}.
\end{remark}



\section{Proof of Theorem \ref{Thm1}.}\label{sectionThm1proof}
We will need the following easy estimate in the proof of Theorem \ref{Thm1}.
\begin{remark}\label{density}
For a given $a>0$ let $S=S_a(\{c_i\})=\{i\in\mathbb{N}|\,
c_i<a\},$ and $s_n=\#(S_a\cap\{i\leq{n}\}).$ If
condition (\ref{product}) holds then
\begin{equation}\frac{s_n}{n}\to1.
\end{equation}
\end{remark}
\begin{proof}
From the usual inequality between
geometric and arithmetic means
$$\sqrt[n]{\prod_{i=1}^n{(1-c_i)}}\leq\frac{1}{n}\sum_{i=1}^n(1-c_i)\leq{1}$$
we get that $\frac{1}{n}\sum_{i=1}^n(1-c_i)\rightarrow{1}$ or,
equivalently, $\frac{1}{n}\sum_{i=1}^n{c_i}\rightarrow0$.
\end{proof}

\begin{remark}\label{RemA}
If we take $c_i\to0$ such that $\sum_ic_i=\infty$ then the corresponding middle interval Cantor set would be an example of a set from Theorem \ref{A}. Indeed, if $c_i\to0$ then $$\frac{1}{n}\sum_{i=1}^n (1-c_i) \asymp \frac{1}{n}\sum_{i=1}^n \log(1-c_i)\to1$$
and then by  (\ref{dim}) $\dim_H E(\{c_i\})=1$. From $\sum_{i}c_i=\infty$ follows that the set has zero measure.  Also, a middle interval Cantor
set $E(\textbf{c})$ is uniformly perfect if and only if there is a constant $C$ such that $c_i<C<1, \forall i\in\mathbb{N}$.

\end{remark}

One of the main tools for proving Theorem \ref{Thm1} will be the
following lemma from ~\cite{H}.

\begin{lemma}
If $f:X\to{Y}$ is $\eta$-quasisymmetric and if
$A\subset{B}\subset{X}$ are such that
$0<\diam{A}\leq\diam{B}<\infty$, then $\diam{f(B)}$ is finite and

\begin{equation}\label{distortion}
\frac{1}{2\eta\left(\frac{\diam{B}}{\diam{A}}\right)}\leq\frac{\diam{f(A)}}{\diam{f(B)}}
\leq\eta\left(\frac{2\diam{A}}{\diam{B}}\right).
\end{equation}
\end{lemma}

By distance between sets below we mean Hausdorff distance: if
$Y,Z\subset{X}$ then
\begin{equation*}
\dist_X(Y,Z)=\inf\{\dist_X(y,z)|\, y\in{Y}, z\in{Z}\}.
\end{equation*}
We will need a different version of (\ref{distortion}).

\begin{lemma}Suppose $X=X_1\cup{X_2}$, with $X_1, X_2$ compact and
$\dist{(X_1,X_2)}>0$. Then
\begin{equation}\label{distortion2}
\frac{1}{2\eta\left(\frac{\diam{X}}{\dist(X_1,X_2)}\right)}
\leq\frac{\dist(f(X_1),f(X_2))}{\diam{f(X)}}
\leq\eta\left(2\frac{\dist(X_1,X_2)}{\diam{X}}\right).
\end{equation}
\end{lemma}

\begin{proof}
Suppose $x_1\in{X_1}$ and $x_2\in{X_2}$ are such that
$\dist{(X_1,X_2)}=d_X(x_1,x_2)$. This is possible since $X_1$ and
$X_2$ are compact. Let $A=\{x_1,x_2\}$ then right hand inequality in
(\ref{distortion}) implies
\begin{eqnarray*}
\frac{\dist(f(X_1),f(X_2))}{\diam{f(X)}}&\leq&\frac{\dist{(f(x_1),f(x_2))}}{\diam{f(X)}}
\leq\eta\left(\frac{2\dist{(x_1,x_2)}}{\diam{X}}\right)\\
&=&\eta\left(2\frac{\dist(X_1,X_2)}{\diam{X}}\right).
\end{eqnarray*}
To obtain the other inequality of (\ref{distortion2}) take
$y_1\in{f(X_1)}, y_2\in{f(X_2)}$ in such a way that
$\dist(f(X_1),f(X_2))=d_Y(y_1,y_2)$. Let $x_i'=f^{-1}(y_i)$. Now
take $A=\{x_1',x_2'\}$. Then again using \ref{distortion} we get
\begin{eqnarray*}
\frac{\dist(f(X_1),f(X_2))}{\diam{f(X)}}\geq\frac{1}{2\eta\left(\frac{\diam{X}}{\dist(x'_1,x'_2)}\right)}
\end{eqnarray*}

Since ${d_X(x_1',x'_2)}\geq\dist(X_1,X_2)$ and since $\eta$ is
increasing we obtain
\begin{equation*}
\frac{1}{2\eta\left(\frac{\diam{X}}{\dist(x'_1,x'_2)}\right)}\geq
\frac{1}{2\eta\left(\frac{\diam{X}}{\dist(X_1,X_2)}\right)}
\end{equation*}
Combining this with the previous inequality gives
(\ref{distortion2}).
    \end{proof}

Theorem \ref{Thm1} follows from the following result and the mass distribution principle
\begin{lemma}
Suppose $f:E\rightarrow Y$  is a power-quasisymmetric
homeomorphism. Then for every $d<1$ there is a measure $\mu$ on $Y$ satisfying
$$\mu(B(y,r))\leq{Cr^d}$$ for some constant $C>0$ all $r>0$ and all
$y\in{Y}$. Constant $C$ does not depend on $y$ and $r$.
\end{lemma}
To simplify the notation below we write $f(E_{n,j})$ for $f(E_{n,j}\cap E)$ (we don't assume that $f$ extends to the real line). We will prove the lemma in several steps.
First we will show that there is a measure $\mu$ on $\bigcap_n\bigcup_{\calE_n} E_{n,j}\subset E$ such that
\begin{equation}\label{growth1}
\mu(f(E_{n,j}))\leq C\diam{f(E_{n,j})},
\end{equation}
for some non-zero finite constant $C$ independent of $n$ and $j$.
\begin{proof}[Proof of \ref{growth1}.]
Every interval $E_{n,j}\in\calE_{n}$
has one ``parent" interval, denoted by $\tilde{E}_{n,j}\in\calE_{n-1}$, containing $E_{n,j}$, and one ``sibling" interval $E_{n,j}'\in\calE_n$ which has the same ``parent".
This notations will
also be used for $f(E_{n,j})$: for an $I\subset Y$ of the form
$I=I_{n,j}=f(E_{n,j})$ we will denote $\tilde I_{n,j}=f(\tilde E_{n,j})$ and $I'_{n,j}=f(E_{n,j}')$.
\subsection{Construction of the measure.}
Now define $\mu$ as follows. Pick $E_0\in\calE_{0}$ and let $$\mu(f(E))=1.$$ For any $I\subset Y$ of
the form $I=f(E_{n,j})$, where $E_{n,j}$ is a ``descendant" of $E_0$ let:
\begin{eqnarray}\label{measure1}
\mu(I)=\frac{\diam^d{I}}{\diam^d{I}+\diam^d{I'}}\mu(\tilde I).
\end{eqnarray}
Given such an interval $I$ there is a unique sequence of nested
subsets
$$I=I_n\subset{I_{n-1}}\subset{I_{n-2}}\subset\ldots\subset{I_2}\subset{I_1}\subset{I_0}={Y}$$
containing it, so that $I_{k-1}=\tilde{I_k}$. By induction we have
\begin{equation*}\label{measure}
\begin{split}
\frac{\mu(I)}{\diam^d{I}}&=\frac{\mu(I_n)}{\diam^d{I_n}}\\
&=\frac{1}{\diam^d{I_n}+\diam^d{I'_n}}\cdot\frac{\diam^d{I_{n-1}}}{\diam^d{I_{n-1}}+\diam^d{I'_{n-1}}}\\
&\quad\quad\cdot\ldots\cdot\frac{\diam^d{I_{1}}}{\diam^d{I_{1}}+\diam^d{I'_{1}}}\mu(I_0).
\end{split}
\end{equation*}

Since $\diam(A\cup{B})\leq\diam{A}+\dist(A,B)+\diam{B}$ we have
\begin{equation}
\begin{split}
\frac{\mu(I)}{\diam^d{I}}&\leq\prod_{i=1}^n\frac{(\diam I_i+\dist
(I_i,I'_i)+\diam{I'_i})^d}{\diam^d{I_i}+\diam^d{I'_i}}.
\end{split}
\end{equation}

Let
\begin{equation}
p_i=\frac{(\diam I_i+\dist
(I_i,I'_i)+\diam{I'_i})^d}{\diam^d{I_i}+\diam^d{I'_i}}\\
\end{equation}
To prove (\ref{growth1}) we need to show that $\prod_{i=1}^n p_i\to0$ as
$n\to\infty$. Indeed, if this is the case then $\exists C<\infty$
s.t. $\prod_{i=1}^n p_i<C, \forall n\in\mathbb{N}$. Now, to prove
$\prod_{i=1}^n p_i\to0$ we will need the following estimates.

\begin{lemma}[Small gaps]\label{smallgaps}
$\,\,\exists \,\,a>0, {C_1}<1$ s.t $c_i<a\Rightarrow{p_i<C_1<1}.$
\end{lemma}
\begin{lemma}[Large gaps]\label{largegaps}
$\,\,\exists \,\,C_2>1$ s.t.  $p_i<\frac{C_2}{(1-c_i)^{d/\alpha}},\forall
i.$
\end{lemma}

Let us prove the theorem assuming these two lemmas. First of all
\begin{equation*}
\begin{split}
\prod_{i=1}^np_i &\leq \prod_{\{i\leq n|
c_i<a\}}C_1\prod_{\{i\leq n| c_i\geq{a}\}}\frac{C_2}{(1-c_i)^{d/\alpha}}\quad\quad(\mbox{by the two lemmas})\\
&\leq \quad
C_1^{s_n}\frac{C_2^{n-s_n}}{{\prod_{i=1}^n(1-c_i)^{d/\alpha}}}\quad\quad(\mbox{where
$s_n$ is like in Corollary \ref{density}}).
\end{split}
\end{equation*}
Now, if $C_1<1$ and $s_n/n\to1$ then for every number $C_2<\infty$
there is a $C_3<1$ and $N\in\mathbb{N}$ s.t. for $n>N$
$$C_1^{s_n}C_2^{n-s_n}\leq{C_3^n}.$$
Hence
\begin{equation*}
\prod_{i=1}^np_i \leq
\left(\frac{C_3}{\sqrt[n]{\prod_{i=1}^n(1-c_i)^{d/\alpha}}}\right)^n.
\end{equation*}
Since $\sqrt[n]{\prod_{i=1}^n(1-c_i)}\to1$ and $C_3<1$ it follows
that $\prod_{i=1}^np_i\to0$.


\subsection{Small gaps.}

\begin{proof} [Proof of lemma \ref{smallgaps}]Recall that for a given $a>0$ we had $$S_a=\{i\in\mathbb{N}|\,
c_i<a\}, S_n=S_a\cap\{i\leq{n}\}, s_n=\textrm{card}(S_n).$$ Without
loss of generality we can assume $a<1/2$.

Suppose now $i\in{S_a}$. We find it easier to estimate $p_i^{-1}$
from below.
\begin{align*}\label{lower bound1}
p_i^{-1}
&=\frac{\diam^d{I_i}+\diam^d{I'_i}}{(\diam{I_i}+\diam{I'_i})^d}\cdot\frac{(\diam{I_i}+\diam{I'_i})^d}{(\diam{I_i}+\dist(I_i,I'_i)+\diam{I'_i})^d}\\
&\geq\frac{\diam^d{I_i}+\diam^d{I'_i}}{(\diam{I_i}+\diam{I'_i})^d}\cdot\left(1-\frac{\dist(I_i,I'_i)}{\diam{I_{i-1}}}\right)^d\\
&\geq{\frac{\diam^d{I_i}+\diam^d{I'_i}}{(\diam{I_i}+\diam{I'_i})^d}}\cdot(1-\eta(2c_i))^d \tag{by (\ref{distortion2})}\\
&=\frac{1+\left(\frac{\diam{I'_i}}{\diam{I_i}}\right)^d}{\left(1+\frac{\diam{I'_i}}{\diam{I_i}}\right)^d}\cdot(1-\eta(2c_i))^d.
\end{align*}
We will show that the the first term in this product is bounded
below by a constant strictly greater than $1$. To do that, first
note that there is a constant $1<D(\eta,M)<\infty$ so that
$D^{-1}<\diam{I_i}/\diam{I'_i}<D.$ Indeed,
\begin{align}
\frac{\diam{I_i}}{\diam{I'_i}}\geq\frac{\diam{I_i}}{\diam{I_{i-1}}}
\geq\frac{1}{2\eta\left(\frac{\diam{E_{i-1}}}{\diam{E_i}}\right)} \tag{by (\ref{distortion2})}
\end{align}
Since $c_{i-1}<1/2$ it follows that $$\dist (E_i,E_i') \leq \diam E_i +\diam E_i'\leq (1+M)\diam E_i.$$
Therefore
$$\frac{\diam{E_{i-1}}}{\diam E_i}\leq \frac{\diam E_i +\dist (E_i,E_i') +\diam E_i'}{\diam E_i}\leq 2(1+M),$$
and hence
$$\frac{\diam{I_i}}{\diam{I'_i}}\geq\frac{1}{2\eta(2(1+M))}>0.$$
The second inequality follows by symmetry.

Considering the function $x\mapsto\frac{1+x^d}{(1+x)^d}$ for $d<1$
one can easily see that on an interval $[D^{-1}, D]$ its smallest
value is attained at $D$ and is strictly larger than $1$. We will
denote this value by $C_4=C_4(\eta,d)>1$. Therefore
\begin{equation}
p_i^{-1} \geq{C_4}(1-\eta(2c_i))^{d}\geq{C_4}(1-\eta(2a))^d.
\end{equation}
Since $\eta$ is increasing and $c_i<a$. Now, $\eta(t)\rightarrow0$
as $t\rightarrow0$. Therefore we can always choose $a$ small enough
so that $C_4(1-\eta(2a))^d>1$. So finally we conclude that there is
an $a$ so that for $i\in S_a$ one has $p_i^{-1}\geq C_5>1$. Equivalently
$p_i$ is bounded from above by a constant strictly less than $1$.\end{proof}
\begin{remark}
Note that we haven't yet used the fact that $f$ is power quasisymmetric.
\end{remark}

\subsection{Large gaps.}

\begin{proof}[Proof of lemma \ref{largegaps}]
Since
$\diam{I_i},\diam{I'_i},\dist(I_i,I_i')<\diam{I_{i-1}}$, we have

\begin{align*}
p_i
&=\frac{(\diam{I_i}+\dist(I_i,I'_i)+\diam{I'_i})^d}{\diam^d{I_i}+\diam^d{I'_i}}\\
&\leq\frac{3^d\diam^d{I_{i-1}}}{\diam^d{I_i}+\diam^d{I'_i}}={3^d}\left[\left(\frac{\diam{I_i}}{\diam{I_{i-1}}}\right)^{d}
+\left(\frac{\diam{I'_i}}{\diam{I_{i-1}}} \right)^{d}\right]^{-1}\\
&\leq{3^d}\frac{\eta^d\left(\frac{\diam E_{i-1}}{\diam E_i}\right)}{2} \tag{by (\ref{distortion2})}.
\end{align*}

From (\ref{ratio}) we have
$$\frac{\diam E_{i}}{\diam E_{i-1}}\geq {1-c_i-\frac{\diam E_i'}{\diam E_{i-1}}}\geq 1-c_i-M\frac{\diam E_i}{\diam E_{i-1}}$$
and therefore
$$\frac{\diam E_{i-1}}{\diam E_{i}}\leq\frac{1+M}{1-c_i}.$$
It follows that
$$p_i\leq \frac{3^d}{2} \eta^d\left(\frac{1+M}{1-c_i}\right).$$
Now, since $\eta(t)\leq{C\max\{t^{1/\alpha},t^{\alpha}\}}$, the last
inequality yields
\begin{equation*}
p_i\leq\frac{C_2(d,\alpha)}{(1-c_i)^{d/\alpha}}.
\end{equation*}\end{proof}
As shown before this completes the proof of (\ref{growth1}).
\end{proof}

\subsection{}To complete the proof of Lemma \ref{growth1} and Theorem \ref{Thm1} we need to show that a similar
estimate holds for any ball $B=B(y,r)$ with $y\in{Y}$. First we show that (\ref{growth1}) implies the following lemma.
\begin{lemma}\label{growth2}
There is a constant $C$ such that
for any interval $J\subset\mathbb{R}$ we have
$$\mu(f(J\cap E))\leq C [\diam f(J\cap E)]^d$$
\end{lemma}
\begin{proof}Note first that for every $J\subset\mathbb R$ there are two (or one) intervals $E_1,E_2\in\bigcup_n\calE_n$ such that
$$E_1, E_2 \subset J \quad \mbox{and} \quad J\cap E\subset \tilde E_1 \cup \tilde E_2.$$
Indeed, consider the collection
$\calE_J=\left\{E_{n,j}\in\bigcup_n\calE_n: E_{n,j}\subseteq J, \,\, \mbox{but} \,\, \tilde E_{n,j}\nsubseteq J\right\},$
in other words, the collection of intervals $E_{n,j}$ which are contained in $J$ with parents $\tilde E_{n,j}$ that are not. Since every interval $E_{n,j}\subset J$ has an ``ancestor'' in $\calE_{J}$ it follows that
$$J\subset \bigcup_{\calE_J}E_{n,j}.$$
Now, choose $E_1\in\calE_J$ so that $\diam\tilde{E_1}\geq\diam \tilde E_{n,j},\, \mbox{for any} \, E_{n,j}\in\calE_J.$
If $\tilde E_1\supset J$ then we are done ($E_2=\emptyset$). If not, consider $\calE_{J\setminus\tilde E_1}$ and choose $E_2$ from this collection in a similar fashion, i.e. $\diam\tilde E_2\geq \diam\tilde E_{n,j}, \mbox{ for any } E_{n,j}\in\calE_{J\setminus\tilde E_1}.$
Since for every $E_{n,j}\in\calE_J$  its parent $\tilde E_{n,j}$ contains at least one of the end points of $J$ it means it intersects either $\tilde E_1$ or $\tilde E_2$ and therefore must be contained in one of them (since every two elements of $\calE$ are either disjoint or one of them contains the other one). Therefore $J\cap E\subset \tilde E_1\cup \tilde E_2$.

Just as before let $E_1'$ and $E_2'$ be the siblings of $E_1$ and $E_2$ respectively.
Note that if $J\cap E_i'=\emptyset$ then $J\cap\tilde{E}_i=E_i$. Therefore we need to consider the contribution of $E_i'$ only if  $J\cap E_i'\neq\emptyset$ in which case, since $\diam E_i'\leq M\diam E_i\leq M\diam J$, we obtain
$$ E_i'\subset 2 M J,$$
where $2MJ$ is just the dilation of $J$ by $2M$. Now from (\ref{growth1}) it follows that
\begin{equation*}
\begin{split}
\mu(f(J\cap E))&=\sum_{i=1,2}\mu(f( E_i))+\mu(f( E_i'))\\
&\leq C\sum_{i=1,2} [\diam f(E_i)]^d+[\diam f(E_i')]^d\\
\end{split}
\end{equation*}
Since
\begin{equation*}
\diam f(E_i),\diam f(E_i')\leq\diam f(2MJ\cap E),
\end{equation*}
and by (\ref{distortion}) we have
\begin{equation*}
\diam f(2MJ\cap E)\leq 2\eta(2M) \diam f(J\cap E)
\end{equation*}
it follows that
$$\mu(f(J\cap E))\leq{C}[\diam f(J\cap E)]^d$$ for some constant $C$ and any interval $J\subset{R}$. \end{proof}

\begin{proof}[Proof of Lemma \ref{growth1}]By quasisymmetry there is a number $1\leq H<\infty$ such that for every $y\in Y$ and $r>0$
$$B(x,R)\subset f^{-1}(B(y,r))\subset B(x,HR),$$
where $x=f^{-1}(y)$.
Therefore
\begin{align*}
\mu(B(y,r))&\leq \mu (f(B(x,HR)))\\
&\leq C[\diam f(B(x,HR))]^d  \tag{by Lemma (\ref{growth1})}\\
&\leq  C[\diam f(B(x,R))]^d \tag{by (\ref{distortion})}\\
&\leq C r^d \tag{$f(B(x,R))\subset B(y,r)$}.
\end{align*}
As we noted before it
follows that $\dim_H(f(E))\geq1$ since $d$ could be chosen as close to
$1$ as one would like.
\end{proof}

\section{Modulus and Conformal dimension}\label{sectionModulus}
As was shown by Tyson in ~\cite{Tyson} one of the main obstructions for lowering the Hausdorff dimension of a space by quasisymmetric maps is the existence of a large family of curves in it. Even though we do not use it below our proof of Theorem \ref{Thm2} is modeled on the proof of Tyson's result given by Bonk and Tyson, see ~\cite{H} Theorem $15.10$.

Recall that measure $\mu$ is said to be \textit{doubling} if there is a number $C$ such that for every ball $B_r$
 $$\mu(B_{2r})\leq C\mu (B_r).$$
A metric measure space $(X,\mu)$ is doubling if $\mu$ is doubling.

\begin{theorem}Suppose $(X,\mu)$ is a doubling metric measure space such that $$\mu(B_r)\lesssim r^d$$ for every ball $B_r\subset X$ of radius $0<r<\diam X$. If there is a curve family $\G$ in $X$ such that $\mod_d \G>0$ then $\confdim X\geq d$.
\end{theorem}

Let us recall that the $d$-modulus of a family of curves $\G$ in $X$ is defined as
$$\mod_d \G = \inf\left\{ \int_X \rho^d d\mu: \, \int_{\g}\rho d s\geq 1, \forall \g\in \G\right\},$$
where $ds$ denotes the arclength element. We refer to ~\cite{H} for further details on modulus of a curve family and the discussion of the theorem of Tyson.

In this section we will give a lower bound on the conformal dimension of a space in terms  of a modulus of a system of measures due to  Fuglede, see ~\cite{Fug}. The need for this comes from the fact that the sets we will be dealing with may have $0$ Hausdorff $1$-measure. In the proof we will need the notion of the discrete modulus of a family of subsets of $X$ which is in essence due to Heinonen and Koskela, see ~\cite{HK}. Below we give the definitions of various moduli formulate the main result, Theorem \ref{Thm2}, and show how Theorem \ref{B} follows from it.

\subsection{Modulus of a system of measures}
Let $(X,\mu)$  be a measure space.
Let $\mathbf{E}$ be a collection of measures on $X$ the domains of which contain the domain of $\mu$. A measurable function $\rho:X\to\mathbb R$ is said to be admissible for the system of measures $\mathbf{E}$ if for every $\la\in\mathbf{E}$ $$\int_E \rho d\la\geq 1.$$
Next we define the $p$-modulus of $\mathbf{E}$ as
$$\mod_p(\mathbf{E})=\inf \int_X\rho^p d\mu,$$
where $\inf$ is taken over all $\mathbf{E}$-admissible functions $\rho$.

Just like the usual modulus of a family of curves the modulus of a system of measures is monotone and sub-additive, see ~\cite{Fug}.

\begin{lemma}
The $p$-modulus is monotone and countably subadditive:
\begin{eqnarray}
\mod_p\mathbf{E}&\leq&\mod_p\mathbf{E}', \mbox{ if } \mathbf{E}\subset\mathbf{E}',\\
\mod_p\mathbf{E}&\leq&\sum_i\mod_p\mathbf{E}_i, \mbox{ if } \mathbf{E}=\bigcup_{i=1}^{\infty}\mathbf{E}_i. \label{subadd}
\end{eqnarray}
\end{lemma}

\subsection{Discrete modulus}
Let $\calE=\{E\}$  be a collection of subsets of $X$.
Let $\calB=\{B\}$ be a cover of  $X$ by balls and $v:\calB\to[0,\infty)$ a function.
The pair $(v,\calB)$ is  \textit{admissible for}  $\calE$ if
$$\frac{1}{5}B \cap\frac{1}{5} B'=\emptyset,$$
whenever $B\neq B'$, and $$\sum_{\frac{1}{5}B\cap E\neq \emptyset} v(B)\geq1$$
for every $E\in\calE$.

For $\delta>0$ set
$$\dmod_{p}^{\d}=\inf \sum_{B\in{\calB}} v(B)^p,$$
where the infimum is over all pairs $(v,\calB)$ which are admissible for $\calE$ and such that $\diam B\leq\d$ for every $B\in\calB$. The \textit{discrete $p$-modulus of} $\calE$ is
$$\dmod_p(\calE)=\lim_{\d\to0} \dmod_{p}^{\d} (\calE).$$

The need for the disjointness property in the definition of admissibility comes from the following covering lemma, see
for instance ~\cite{Mattila} Theorem $2.3$.

\begin{lemma}[\textbf{Covering Lemma}]\label{covlemma}
Every family $\calB$ of balls of bounded diameter in a compact metric space $X$
contains a countable subfamily of disjoint balls $B_i\subset\calB$ such that
$$\bigcup_{B\in\calB}B\subset\bigcup_{i}5B_i.$$
\end{lemma}

\begin{remark}
Even though the monotonicity of the discrete modulus is easy to see we do not know if the analogue of (\ref{subadd}) is true in this case.
\end{remark}

\subsection{Conformal dimension and Fuglede modulus}
\begin{theorem}\label{Thm2}
Let $p>q>1$ and $(X,\mu)$ be a doubling metric measure space. Suppose there is a constant $0<C<\infty$ such that for every ball $B_r\subset X$
\begin{equation}\label{cond1}
\mu(B_r)\leq C r^{p}.
\end{equation}
Let $\calE$ be a collection of subsets of $X$ such that
\begin{align}
\confdim E &\geq1, \forall E\in\calE, \label{cond2}
\end{align}
If there is  a system of measures $\mathbf{E}=\{\la_E\}$ associated to $\calE$ so that $$\mod_q {\mathbf{E}}>0$$ and for every $s>1$ there are constants $C_1=C_1(s)$ and $C_2=C_2(s)$ such that $\forall E\in\mathcal{E}$ and $B_r\subset X$
\begin{equation}\label{cond4}
\lambda_E(B_r\cap E)\geq C_1 r^{s},
\end{equation}
provided $\frac{1}{C_2}B_r\cap E\neq\emptyset,$ then
$$\confdim X\geq q.$$
\end{theorem}

The proof of the theorem is given in the next section. Here we show how Theorem \ref{B} follows from Theorem \ref{Thm2}.

\begin{corollary}
If $E\subset\mathbb{R}$ is a set of conformal dimension $1$ which supports a measure $\la_E$ such that for every $\eps>0$ there is a constant $C$ so that
 $$\frac{1}{C}R^{1+\eps}\leq \la_E(B_R) \leq C R^{1+\eps}$$
 then for every Borel set $Y\subset\mathbb{R}^n$
 $$\confdim(E\times{Y})\geq\dim_H E\times{Y}.$$
\end{corollary}
\begin{proof}
Let $d<\dim_H Y$. By Frostman's  lemma for every $\eps$ such that $2\eps\in(0, \dim_H Y - d)$ there is a measure $\nu$ on $Y$ such that $\nu(Y)>0$ and $\nu(B_R)\lesssim R^{d+2\eps}$ for every ball $B_R\subset Y.$ Let $\mu=\la_E \times \nu$.
 Then there is a constant $0<C<\infty$ such that  $$\mu(B_R)\leq C R^{1-\eps} R^{d+2\eps}=R^{1+d+\eps}$$ for every $B_R\subset E\times Y$.

Let $\calE=\{E\times\{y\}: y\in Y\}$ and $\la_{E\times\{y\}}(U\times\{y\})=\la_E(U)$ for every $\la_E$ measurable $U\subset E$. Define
$$\mathbf{E}=\{ \la_{E\times\{y\}}: y\in Y\}.$$
The proof would be complete if we could show that $\mod_{1+d} \mathbf{E} >0$. Indeed, Theorem \ref{Thm2} would imply then that $\confdim (E\times Y)\geq 1+d$ for every $d<\dim_H Y$ and therefore $\confdim (E\times Y)\geq 1+\dim_H {Y}$.

The argument for $\mod_{1+d} \mathbf{E}>0$ is standard and we include it only for completeness.
Take $\rho:E\times{Y}\to\mathbb{R}_+$ s.t. $\int_E \rho(x,y) d\la_E\geq1, \forall y\in{Y}$. By H\"older's inequality we get that
$$\int_{E}\rho^{1+d}(x,y)dx\geq1, \forall y\in{Y}.$$
Integrating both sides of the inequality with respect to $\nu$ we obtain
$$\int_{E\times Y}\rho^{1+d}(x,y) d\la_E\times d\nu\geq \nu(Y)>0.$$
Therefore $\mod_{1+d}\mathbf{E}\geq\nu(Y)>0$.
\end{proof}

\begin{remark}
It is not hard to see that uniformly perfect middle interval Cantor sets  satisfy the conditions of the previous corollary. In fact the measure which gives equal mass to every interval of the same length is an example of a measure which satisfies the required inequalities.
\end{remark}

\section{Proof of Theorem \ref{Thm2}}\label{sectionThm2}
Theorem \ref{Thm2} would follow from the following two lemmata.

\begin{lemma}\label{dmod=0} Let $t>0$ and suppose  $\calE$ is a collection of subsets in $X$ such that for some $\delta>0$
$$H_t^{\delta} E\geq c>0, \forall E \in \calE.$$  Then $\dmod_q \calE=0$ for every $q>\frac{1}{t}\dim_H {X}$.
\end{lemma}
\begin{lemma}\label{moduli inequality} If conditions (\ref{cond1}) and (\ref{cond4}) of Theorem \ref{Thm2} are satisfied. Then for every $q<p$ there is a constant $C<\infty$ such that
$$\mod_q{\mathbf E}\leq C \dmod_q f(\calE).$$
\end{lemma}
Before proving the lemmas let us prove the theorem assuming they are true.

\begin{proof}[Proof of Theorem \ref{Thm2}] Suppose $f$ is a quasisymmetric map such that $\dim_H f(X)<q$.
Choose $t<1$ so that $qt>\dim_H X$. Since $\dim_H f(E)\geq 1$ for every $E\in\calE$ it follows that $H_t^{1/k}f(E)\to\infty$ as $k\to\infty$. Let $k_E\in\mathbb N$ denote the smallest integer such that $$H_t^{1/k_E}f(E)\geq1.$$
Let $$\calE_j=\{E\in\calE: k_E\geq j\}$$
and
$$\mathbf E_j=\{\la_E\in\mathbf{E}: E\in\calE_j\}.$$
Then $$\mathbf E=\bigcup_{j=1}^{\infty}\mathbf E_j$$
and therefore
\begin{align*}
\mod_q\mathbf E
&\leq\sum_{i=1}^{\infty} \mod_q\mathbf E_i \tag{by (\ref{subadd})}\\
&\leq C\sum_{i=1}^{\infty} \dmod_q f(\calE_i), \tag{by Lemma \ref{moduli inequality}}
\end{align*}
where $f(\calE_i)$ is the image of the family $\calE_i$. Lemma \ref{dmod=0} implies that $\dmod_q f(\calE_j)=0$ and it follows that $\mod_q \mathbf{E}=0$ which contradicts our assumption.
\end{proof}

\begin{remark}
Here are some questions which naturally arise and would simplify and generalize the proof above.
\begin{itemize}
\item[1.] {Is $\dmod_q$ countably subadditive?}
\item[2.] {Is $\dmod_q$ a quasisymmetric quasi-invariant?}
\item[3.] {Is $\dmod_q(E\times Y)>0$ for $q=\dim_H(E\times Y)$ (or under what conditions this is true)?}
\end{itemize}
\end{remark}

\subsection{Proof of Lemma \ref{dmod=0}}

\begin{proof}
Let  $q'$ be such that $\dim_H X<q'<tq$. For every $\delta>\eps>0$ there is a covering $\calB$ of $X$ by balls $B_1, B_2,\ldots$ with radii  $r_1,r_2,\ldots$ such that
$\frac{1}{5}B_i\cap \frac{1}{5} B_j =\emptyset$, for $i\neq j$ and $$\sum_{i} r_i^{q'}<\eps.$$
Let $v(B_i)=r_i^t.$
Since $r_i<\delta$ it follows that for every $E\in\calE$ we have
$$\sum_{B\cap E\neq\emptyset}v(B)\geq H^{\delta}_t(E)\geq1$$
and so $(v,\calB)$ is admissible for $\calE$. Now,
$$\sum_{i} v(B_i)^q = \sum_{i} r_i^{tq} \leq \sum_i r_i^{q'}<\eps.$$
Therefore $\dmod_q\calE<\eps$.
\end{proof}

\subsection{Proof of Lemma \ref{moduli inequality}} Below we will need the following well known inequality, see ~\cite{H} or ~\cite{Bojarski} Lemma $4.2$ in the case of $\mathbb{R}^n$,  which is a consequence of the boundedness of the Hardy-Littlewood maximal operator.

\begin{lemma}Suppose $\mathcal{B}=\{B_i,B_2,\ldots\}$ is a countable collection of balls in a doubling metric measure space $(X,\mu)$ and $a_i\geq0$ are real numbers. Then there is a positive constant $C$ such that
\begin{equation}
\int_X \left(\sum_{\mathcal{B}}a_i\chi_{AB_i}(x)\right)^p d\mu \leq C(A,p,\mu) \int_X \left(\sum_{\mathcal{B}}a_i\chi_{B_i}(x)\right)^p d\mu \tag{$\star$}
\end{equation}
for every $1<p<\infty$ and $A>1$.
\end{lemma}

\begin{proof}[{Proof of Lemma \ref{moduli inequality}}]
It is clear that all we need to show is $$\mod_q\mathbf{E}\leq C\dmod_q^{\delta}f(\calE)$$ for some $\delta\in(0,1)$. For that suppose $(v,\calB')$ is an $f(\calE)$-admissible pair, where $\calB'=\{B_i'\}_{i=1}^{\infty}$ is a cover of $f(X)$ by balls  $B_i'$ of radii $r_i'<\delta$ .  Choose $B_i\subset X$ with radius $r_i$ so that
$$\frac{1}{H}B_i\subset f^{-1}\left(\frac{1}{5}B_i'\right)\subset f^{-1}(B_i')\subset B_i,$$ where $H$ is constant depending on $f$ (there is such a constant since $f$ is quasisymmetric). Note that since $\calB'$ is admissible it follows that $\frac{1}{H}B_i\cap\frac{1}{H}B_j=\emptyset$ whenever $i\neq j$.

We want to construct an $\mathbf{E}$-admissible function $\rho$ such that $$\int_X \rho^q d\mu\leq C\sum_{B\in\calB'}v(B')^q.$$
Define
$$\rho(x)=\sum_i\frac{v(B_i')}{[\diam B_i]^{s}}\chi_{C_2B_i}(x),$$
where $C_2$ is as in the formulation of Theorem \ref{Thm2}. Then for every $E\in\mathcal{E}$  the following holds
\begin{align*}
\int_E\rho  d\lambda_E &=\int_E\sum_i\frac{v(B_i')}{[\diam B_i]^{s}}\chi_{C_2B_i}(x) d\la_E
\geq \int_E\sum_{i: \,B_i\cap E\neq\emptyset} \frac{v(B_i')}{[\diam{B_i}]^{s}}\chi_{C_2B_i}(x) d\la_E\\
& = \sum_{i: \,B_i\cap E\neq\emptyset} \frac{v(B_i')}{[\diam{B_i}]^{s}}\int_{E\cap C_2B_i}d \lambda_E =\sum_{i: \,B_i\cap E\neq\emptyset} \frac{v(B_i')}{[\diam{B_i}]^s} \lambda_E ({E\cap C_2B_i})\\
&\geq \frac{1}{C_1} \sum _{i:f(E)\cap B_i'\neq\emptyset}{v(B_i')} \geq \frac{1}{C_1}  \nonumber
\end{align*}
It follows that $$\mod_{q}(\mathbf{E})\leq C_1^q\int_X\rho^q d\mu.$$

Next, take $s>1$ so that $qs<p$. Then we have
\begin{align}
\int_X \rho^q d\mu &= \int_X \left(\sum_i\frac{v(B_i')}{[\diam B_i]^{s}}\chi_{5B_i}(x)\right)^q d\mu\nonumber \\
&\leq C(5H, q, \mu)\int_X \left(\sum_i \frac{v(B_i')}{[\diam {B_i}]^{s}}\chi_{\frac{1}{H}B_i}(x) \right)^q d\mu \tag{by ($\star$)} \\
& =  C(5H, q, \mu)\sum_i \left(\frac{v(B_i')}{[\diam{B_i}]^{s}}\right)^q \mu\left(\frac{1}{H}B_i\right) \\
&\lesssim \sum_i v(B_i')^q r_i^{p-qs} \tag{by (\ref{cond1})}\\
&\lesssim \sum_iv(B_i')^{q} \tag{ $r_i<\delta<1$}.
\end{align}
Taking infimum over all $f(\calE)$-admissible pairs $(v',\calB')$ we obtain $\mod_q\mathbf{E}\leq C\dmod_q^{\delta}f(\calE)$ for some $C$ independent of $\delta$ and hence $$\mod_q\mathbf{E}\leq C\dmod_qf(\calE)$$ therefore completing the proof.
\end{proof}

%
%
%

\end{document}